\documentclass[11pt,twoside]{article}
\usepackage[authoryear,semicolon]{natbib}
\usepackage[T1]{fontenc}
\usepackage[utf8]{inputenc}
\usepackage{textgreek}
\usepackage{fullpage}
\usepackage[color=yellow]{todonotes}
\usepackage{xifthen}
\usepackage{amsmath,amssymb,amsthm,mathtools,stmaryrd}
\usepackage{proof,mathpartir}
\usepackage{colonequals}
\usepackage{verbatim}
\usepackage{comment}
\usepackage{textcomp}
\usepackage[us]{optional}
\usepackage{color}
\usepackage{url}
\usepackage{graphics}
\usepackage{import}
\usepackage{stackengine}
\usepackage{scalerel}

\AtBeginDocument{\mathcode`\;=\numexpr\mathcode`\;-\string"4000\relax}

\theoremstyle{plain}
\newtheorem{theorem}{Theorem}

\newtheorem{lemma}[theorem]{Lemma}

\newtheorem{exercise}{Exercise}

\theoremstyle{remark}

\newcommand{\eqdef}{\mathrel{\triangleq}}
\newcommand{\isdef}{\eqdef}
\newcommand{\bnfdef}{\coloncolonequals}
\newcommand{\bnfalt}{\mathrel{\mid}}

\newcommand{\parens}[1]{(#1)}
\newcommand{\braces}[1]{\{#1\}}
\newcommand{\sqbracks}[1]{[#1]}

\newcommand{\subst}[3]{[{#1}/{#2}]{#3}}

\newcommand{\const}[1]{\textbf{#1}}

\newcommand{\LF}[1][]{\ensuremath{\mathsf{LF}_{#1}}}
\newcommand{\sortclass}{\ensuremath{\textsf{Sort}}}
\newcommand{\eqclass}[3]{{#2}=_{#1}{#3}}
\newcommand{\piclass}[3]{\braces{{#2}\mathbin{:}{#1}}\,{#3}}
\newcommand{\pisort}[3]{\piclass{#1}{#2}{#3}}
\newcommand{\arrclass}[2]{{#1}\to{#2}}
\newcommand{\lvlclass}{\ensuremath{\mathsf{Lvl}}}
\newcommand{\lamobj}[3]{\sqbracks{{#2}\mathbin{:}{#1}}\,{#3}}

\newcommand{\appobj}[2]{{#1}\,{#2}}
\newcommand{\selfobj}{\bullet}
\newcommand{\zeroobj}{\mathbf{0}}
\newcommand{\succobj}[1]{{#1}\mathbin{\mathbf{+}}\mathbf{1}}
\newcommand{\empctx}{\varepsilon}
\newcommand{\snocctx}[3]{{#1}\mathbin{,}{#2}{:}{#3}}
\newcommand{\appctx}[2]{{#1}\,{#2}}

\newcommand{\issig}[1]{{#1}\;\mathsf{sig}}
\newcommand{\isctx}[2][]{{#2}\;\mathsf{ctx}_{#1}}

\newcommand{\iscls}[3][]{{#2}\vdash_{#1}{#3}\;\mathsf{cls}}
\newcommand{\eqcls}[4][]{{#2}\vdash_{#1}{#3}={#4}\;\mathsf{cls}}
\newcommand{\isobj}[4][]{{#2}\vdash_{#1}{#3}:{#4}}
\newcommand{\eqobj}[5][]{{#2}\vdash_{#1}{#3}={#4}:{#5}}

\allowdisplaybreaks[1]       

\begin{document}

\title{An Equational Logical Framework for Type Theories}
\author{Robert Harper}
\date{\today}

\maketitle{}

\section{Introduction}

A \emph{logical framework} is language for defining logical systems, in particular type
theories.  The definition of a logical system consists of a collection of
\emph{generators} of the \emph{objects} that populate a collection of classifying
\emph{sorts}, and a collection of \emph{equations} that identify objects of those sorts.
The generators specify the sorts and objects that constitute the logical system---say, the
sort of types and, for each type, the sort of its elements, and for each sort of element,
the objects that are its members.  The sorts and objects include the syntactic entities of
the logical system---say, the types and their elements---and are specified using
\emph{higher-order abstract syntax} to express the binding and scopes of
variables~\citep{harper-etal:lf}.  They also include the deductive apparatus of a
logic---its \emph{basic}, \emph{hypothetical}, and \emph{general judgments}, and the
\emph{evidence} for them~\citep{martin1987truth}---and these are also specified using
higher-order representations~\citep{harper-etal:lf}.

The native equality of the framework is a congruence---an equivalence relation compatible
with the generators---and defines substitution of objects for variables in another object.
The defining equations of a logical system enrich the native equality to specify the
behavior of the represented objects---such as the inversion and unicity properties of type
constructors or connectives.  The integration of the defining equations with the native
equations of the framework constrains the \emph{meaning} of the defined system within the
framework in the sense that any interpretation must obey the specified laws.  There being
no limitations on the nature of these equations, the enriched equality judgment of the
framework may or may not be (feasibly or infeasibly) decidable.  A \emph{syntactic logical
  framework}~\citep{harper-etal:lf} is one that presents a logical system using only
generators, and no relations, so that the induced equational theory is the native one,
which is decidable.  A \emph{semantic logical
  framework}~\citep{smith1990programming,uemura19} admits the specification of an
equational theory that may not be (feasibly) decidable.  Each logical system is a problem
of its own, and it is in general difficult to transfer results from one case to another.

This note defines a semantic logical framework suitable for defining a broad---but by no
means comprehensive---class of logical systems, including full-scale dependent type
theories.  It is a dependently typed language with a single Russellian universe of sorts,
an extensional equality types governing the objects of a sort, and closed under the
formation of dependent function types.  The definition of a logical system is a form of
context, called a \emph{signature}, that specifies generators that populate the sorts and
the equality types that govern them.  The \emph{adequacy} of a signature expresses the
intended correspondence between the components of the represented logical system and their
counterpart objects in the logical framework.

\emph{Acknowledgement}  I am grateful to Jon Sterling and Carlo Angiuli for many
discussions about this paper, and to James Andrews and Anja Petkovi\'{c} for corrections.

\section{A Logical Framework}

The syntax of the logical framework \LF{} is given in Figure~\ref{fig:lf-syntax}.  It is a
dependently typed $\lambda$-calculus with the structure specified in the introduction.
The types of \LF{} are called \emph{classes}, $K$ or $S$, and their elements are called
\emph{objects}, $O$ or $S$.\footnote{The double role of the meta-variable $S$ will be
  explained shortly.}  The notation is inspired by \textsf{AUTOMATH}, using square
brackets for $\lambda$-abstraction, curly braces for $\Pi$-types, and juxtaposition for
application.  The binding and scopes of identifiers are understood; all classes and
objects are identified up to renaming of bound variables.  Substitution of an object for a
variable within a class is defined in the usual way up to such renamings.

\begin{figure}[tp]
  
  \begin{displaymath}
    \begin{array}{lrcl}
      \textit{Variables} & X,x \\
      \textit{Classes} & K, S & \bnfdef{} & S \bnfalt \sortclass \bnfalt \piclass{S_{1}}{X}{K_{2}}
                                          \bnfalt \eqclass{S}{O_{1}}{O_{2}} \\
      \textit{Objects} & O, S & \bnfdef & X \bnfalt \selfobj{}
                                          \bnfalt \pisort{S_{1}}{X}{S_{2}} \bnfalt \lamobj{S_{1}}{X}{O_{2}} \bnfalt
                                          \appobj{O}{O_{1}} \\
      \textit{Contexts} & \Gamma & \bnfdef & \empctx{} \bnfalt \snocctx{\Gamma}{X}{K}
    \end{array}
  \end{displaymath}

  \caption{Abstract Syntax of \LF{}}
  \label{fig:lf-syntax}
\end{figure}

\LF{} uses dependent function classes to define the \emph{hypothetico-general judgment
  form} that is central to the definition of many logical systems.  It uses (extensional)
equality classes as a convenient way to present equational theories, and it uses the
objects of class \sortclass{} for the syntactic categories of a logical system.  It is
itself defined in the conventional manner in terms of the hypothetico-general judgment
forms given in Figure~\ref{fig:lf-judgments}.

\begin{figure}[tp]
  \centering
  \begin{tabular}{l@{\qquad}l}
    $\isctx{\Gamma}$ & $\Gamma$ is a context \\[1ex]
    $\iscls{\Gamma}{K}$ & $K$ is a class in context $\Gamma$ \\
    $\isobj{\Gamma}{O}{K}$ & $O$ is an object of class $K$ in context $\Gamma$ \\[1ex]
    $\eqcls{\Gamma}{K}{K'}$ & $K$ and $K'$ are equal classes in context $\Gamma$ \\
    $\eqobj{\Gamma}{O}{O'}{K}$ & $O$ and $O'$ are equal objects of class $K$ in context $\Gamma$
  \end{tabular}
  \caption{Judgment Forms of \LF{}}
  \label{fig:lf-judgments}
\end{figure}

A \emph{signature} $\Sigma$ is a context.  The variables declared in a signature are
written $C$ and $c$ to suggest their role as constants, or generators.  The specialization
of the judgment forms of \LF{} to a particular signature $\Sigma$ are defined in
Figure~\ref{fig:lf-specialized} using concatenation of contexts, which is defined in the
evident manner.

\begin{figure}[tp]
  \centering
  \begin{tabular}{l@{\qquad\qquad}l}
    $\issig{\Sigma}$             & $\isctx{\Sigma}$ \\[1ex]

    $\isctx[\Sigma]{\Gamma}$           & $\isctx{\appctx{\Sigma}{\Gamma}}$  \\[1ex]
    
    $\iscls[\Sigma]{\Gamma}{K}$           & $\iscls{\appctx{\Sigma}{\Gamma}}{K}$ \\
    $\isobj[\Sigma]{\Gamma}{O}{K}$     & $\isobj{\appctx{\Sigma}{\Gamma}}{O}{K}$ \\[1ex]

    $\eqcls[\Sigma]{\Gamma}{K}{K'}$   & $\eqcls{\appctx{\Sigma}{\Gamma}}{K}{K'}$ \\
    $\eqobj[\Sigma]{\Gamma}{O}{O'}{K}$ & $\eqobj{\appctx{\Sigma}{\Gamma}}{O}{O'}{K}$
  \end{tabular}
  \caption{Specialized Judgment Forms of \LF{}}
  \label{fig:lf-specialized}
\end{figure}

The rules defining the \LF{} judgment forms are given in Figures~\ref{fig:lf-is},
\ref{fig:lf-str}, and~\ref{fig:lf-eq} which are to be understood as constituting one
simultaneous inductive definition.  Rule~\textsc{incl} specifies that every sort is itself
a class; the class $\sortclass$ is thus formulated as a Russellian universe of ``small''
classes.  Rule~\textsc{pi}, and associated rules~\textsc{lam} and~\textsf{app}, are
restricted to require that $S_{1}$ be a sort, rather than a general class,\footnote{But
  variables are nevertheless permitted to range over classes.}  and, similarly, equality
classes are limited to equations between objects of a sort.  The class \sortclass{} is
required to be closed under dependent function sorts.  Rule~\textsc{eq} defines the class
of equations between objects of a sort; rule~\textsc{self} specifies that every object is
equal to itself.  Rules~\textsc{obj-cls} and~\textsc{obj-eq-cls} specifies that equal
classes classify the same objects.  Rules~\textsc{app-lam} and~\textsc{lam-app} specify
the inversion principles governing abstraction and application.  Rule~\textsc{reflection}
specifies that the equality class internalizes equality, and rule~\textsc{unicity}
specifies that an equation is ``at most true,'' there being no distinction between two
objects witnessing its truth.

\begin{lemma}[Presuppositions]
  \label{lemma:presup}
  \begin{enumerate}
  \item If\/ $\iscls{\Gamma}{K}$, then $\isctx{\Gamma}$, and if\/ $\eqcls{\Gamma}{K}{K'}$, then $\iscls{\Gamma}{K}$
    and $\iscls{\Gamma}{K'}$.
  \item If\/ $\isobj{\Gamma}{O}{K}$, then $\iscls{\Gamma}{K}$, and if\/ $\eqobj{\Gamma}{O}{O'}{K}$, then
    $\isobj{\Gamma}{O}{K}$ and $\isobj{\Gamma}{O}{K'}$.
  \end{enumerate}
\end{lemma}

\begin{lemma}[Weakening]
  \label{lemma:weak}
  Suppose that $\isctx[\Gamma_1]{\Gamma_2}$.
  \begin{enumerate}
  \item If\/ $\iscls{\Gamma_1}{K}$, then $\iscls{\appctx{\Gamma_{1}}{\Gamma_2}}{K}$, and if\/
    $\eqcls{\Gamma_1}{K}{K'}$, then $\eqcls{\appctx{\Gamma_1}{\Gamma_2}}{K}{K'}$.
  \item If\/ $\isobj{\Gamma_1}{O}{K}$, then $\isobj{\appctx{\Gamma_1}{\Gamma_2}}{O}{K}$, and if\/
    $\eqobj{\Gamma_1}{O}{O'}{K}$, then $\eqobj{\appctx{\Gamma_1}{\Gamma_2}}{O}{O'}{K}$.
  \end{enumerate}
\end{lemma}

\begin{lemma}[Substitution]
  \label{lemma:subst}
  Let $\Gamma\isdef{}\appctx{\appctx{\Gamma_{1}}{X{:}K_{1}}}{\Gamma_{2}}$, and suppose that
  $\isctx{\Gamma}$ and $\isobj{\Gamma_{1}}{O_{1}}{K_{1}}$.
  \begin{enumerate}
  \item If\/ $\iscls{\Gamma}{K_{2}}$, then
    $\iscls{\appctx{\Gamma_1}{\subst{O_1}{X}{\Gamma_2}}}{\subst{O_{1}}{X}{K_{2}}}$, and similarly for class equality.
  \item If\/ $\isobj{\Gamma}{O_{2}}{K_{2}}$, then
    $\isobj{\appctx{\Gamma_1}{\subst{O_1}{X}{\Gamma_2}}}{\subst{O_{1}}{X}{O_{2}}}{\subst{O_{1}}{X}{K_{2}}}$,
    and similarly for object equality.
  \end{enumerate}
\end{lemma}

\begin{lemma}[Functionality]
  \label{lemma:func}
    Let $\Gamma\isdef{}\appctx{\appctx{\Gamma_{1}}{X{:}K_{1}}}{\Gamma_{2}}$, and suppose that
    $\isctx{\Gamma}$ and $\eqobj{\Gamma_{1}}{O_{1}}{O_{1}'}{K_{1}}$.
  \begin{enumerate}
  \item If\/ $\iscls{\Gamma}{K_{2}}$, then
    $\eqcls{\appctx{\Gamma_1}{\subst{O_1}{X}{\Gamma_2}}}{\subst{O_{1}}{X}{K_{2}}}{\subst{O_{1}'}{X}{K_{2}}}$.
  \item If\/ $\isobj{\Gamma}{O_{2}}{K_{2}}$, then
    $\eqobj{\appctx{\Gamma_1}{\subst{O_1}{X}{\Gamma_2}}}{\subst{O_{1}}{X}{O_{2}}}{\subst{O_{1}'}{X}{O_{2}}}{\subst{O_{1}}{X}{K_{2}}}$,
  \end{enumerate}
\end{lemma}

\begin{figure}[tp]
  
  \begin{mathpar}
    
    \inferrule[sort]
    {\isctx{\Gamma}}
    {\iscls{\Gamma}{\sortclass}}

    \inferrule[incl]
    {\isobj{\Gamma}{S}{\sortclass}}
    {\iscls{\Gamma}{S}}

    \inferrule[pi-cls]
    {\isobj{\Gamma}{S_{1}}{\sortclass} \\ \iscls{\snocctx{\Gamma}{X}{S_{1}}}{K_{2}}}
    {\iscls{\Gamma}{\piclass{S_{1}}{X}{K_{2}}}}

    \inferrule[eq-cls]
    {\isobj{\Gamma}{S}{\sortclass} \\ \isobj{\Gamma}{O_{1}}{S} \\ \isobj{\Gamma}{O_{2}}{S}}
    {\iscls{\Gamma}{\eqclass{S}{O_{1}}{O_{2}}}}

  \end{mathpar}

  \begin{mathpar}

    \inferrule[pi-sort]
    {\isobj{\Gamma}{S_{1}}{\sortclass} \\ \isobj{\snocctx{\Gamma}{X}{S_{1}}}{S_{2}}{\sortclass}}
    {\isobj{\Gamma}{\pisort{S_{1}}{X}{S_{2}}}{\sortclass}}

    \inferrule[lam]
    {\isobj{\Gamma}{S_{1}}{\sortclass} \\ \isobj{\snocctx{\Gamma}{X}{S_{1}}}{O_{2}}{K_{2}}}
    {\isobj{\Gamma}{\lamobj{S_{1}}{X}{O_{2}}}{\piclass{S_{1}}{X}{K_{2}}}}

    \inferrule[app]
    {\isobj{\Gamma}{O}{\piclass{S_{1}}{X}{K_{2}}} \\ \isobj{\Gamma}{O_{1}}{S_{1}}}
    {\isobj{\Gamma}{\appobj{O}{O_{1}}}{\subst{O_{1}}{X}{K_{2}}}}

    \inferrule[self]
    {\isobj{\Gamma}{O}{S}}
    {\isobj{\Gamma}{\selfobj}{\eqclass{S}{O}{O}}}

  \end{mathpar}

  \caption{Formation Judgments}
  \label{fig:lf-is}
\end{figure}

\begin{figure}[tp]
  
  \begin{mathpar}
    
    \inferrule[emp]
    {\strut}
    {\isctx{\empctx}}

    \inferrule[snoc]
    {\iscls{\Gamma}{K}}
    {\isctx{\snocctx{\Gamma}{X}{K}}}

    \inferrule[var]
    {\isctx{\appctx{\appctx{\Gamma_{1}}{X{:}K}}{\Gamma_{2}}}}
    {\isobj{\appctx{\appctx{\Gamma_{1}}{X{:}K}}{\Gamma_{2}}}{X}{K}}

  \end{mathpar}

  \begin{mathpar}
    
    \inferrule[cls-rfl]
    {\iscls{\Gamma}{K}}
    {\eqcls{\Gamma}{K}{K}}

    \inferrule[cls-st]
    {\eqcls{\Gamma}{K}{K'} \\ \eqcls{\Gamma}{K''}{K'}}
    {\eqcls{\Gamma}{K}{K''}}

  \end{mathpar}

  \begin{mathpar}

    \inferrule[obj-rfl]
    {\isobj{\Gamma}{O}{K}}
    {\eqobj{\Gamma}{O}{O}{K}}

    \inferrule[obj-st]
    {\eqobj{\Gamma}{O}{O'}{K} \\ \eqobj{\Gamma}{O''}{O'}{K}}
    {\eqobj{\Gamma}{O}{O''}{K}}

  \end{mathpar}

  \begin{mathpar}
    
    \inferrule[obj-cls]
    {\isobj{\Gamma}{O}{K} \\ \eqcls{\Gamma}{K}{K'}}
    {\isobj{\Gamma}{O}{K'}}

    \inferrule[obj-eq-cls]
    {\eqobj{\Gamma}{O}{O'}{K} \\ \eqcls{\Gamma}{K}{K'}}
    {\eqobj{\Gamma}{O}{O'}{K'}}

  \end{mathpar}

  \caption{Structural Judgments}
  \label{fig:lf-str}
\end{figure}

\begin{figure}[tp]
  
  \begin{mathpar}
    
    \inferrule[incl-eq]
    {\eqobj{\Gamma}{S}{S'}{\sortclass}}
    {\eqcls{\Gamma}{S}{S'}}

    \inferrule[pi-class-eq]
    {\eqobj{\Gamma}{S_{1}}{S_{1}'}{\sortclass} \\ \eqcls{\snocctx{\Gamma}{X}{S_{1}}}{K_{2}}{K_{2}'}}
    {\eqcls{\Gamma}{\piclass{S_{1}}{X}{K_{2}}}{\piclass{S_{1}'}{X}{K_{2}'}}}

    \inferrule[eq-class-eq]
    {\eqobj{\Gamma}{S}{S'}{\sortclass} \\ \eqobj{\Gamma}{O_{1}}{O_{1}'}{S} \\ \eqobj{\Gamma}{O_{2}}{O_{2}'}{S'}}
    {\eqcls{\Gamma}{\eqclass{S}{O_{1}}{O_{2}}}{\eqclass{S'}{O_{1}'}{O_{2}'}}}

  \end{mathpar}

  \begin{mathpar}
    
    \inferrule[pi-sort-eq]
    {\eqobj{\Gamma}{S_1}{S_1'}{\sortclass} \\ \eqobj{\snocctx{\Gamma}{X}{S_1}}{S_2}{S_2'}{\sortclass}}
    {\eqobj{\Gamma}{\pisort{S_1}{X}{S_2}}{\pisort{S_1'}{X}{S_2'}}{\sortclass}}

    \inferrule[lam-eq]
    {\eqcls{\Gamma}{S_{1}}{S_{1}'} \\ \eqobj{\snocctx{\Gamma}{X}{S_{1}}}{O_{2}}{O_{2}'}{K_{2}}}
    {\eqobj{\Gamma}{\lamobj{S_{1}}{X}{O_{2}}}{\lamobj{S_{1}'}{X}{O_{2}'}}{\piclass{S_{1}}{X}{K_{2}}}}

    \inferrule[app-eq]
    {\eqobj{\Gamma}{O}{O'}{\piclass{S_{1}}{X}{K_{2}}} \\ \eqobj{\Gamma}{O_{1}}{O_{1}'}{S_{1}}}
    {\eqobj{\Gamma}{\appobj{O}{O_{1}}}{\appobj{O'}{O_{1}'}}{\subst{O_{1}}{X}{K_{2}}}}

  \end{mathpar}

  \begin{mathpar}
    
    \inferrule[app-lam]
    {\isobj{\snocctx{\Gamma}{X}{S_{1}}}{O_{2}}{K_{2}} \\ \isobj{\Gamma}{O_{1}}{S_{1}}}
    {\eqobj{\Gamma}{\appobj{(\lamobj{S_{1}}{X}{O_{2}})}{O_{1}}}{\subst{O_{1}}{X}{O_{2}}}{\subst{O_{1}}{X}{K_{2}}}}

    \inferrule[lam-app]
    {\isobj{\Gamma}{O}{\piclass{S_{1}}{X}{K_{2}}}}
    {\eqobj{\Gamma}{O}{\lamobj{S_{1}}{X}{(\appobj{O}{X})}}{\piclass{S_{1}}{X}{K_{2}}}}

  \end{mathpar}

  \begin{mathpar}
    
    \inferrule[reflection]
    {\isobj{\Gamma}{O}{\eqclass{S}{O_{1}}{O_{2}}}}
    {\eqobj{\Gamma}{O_{1}}{O_{2}}{S}}

    \inferrule[unicity]
    {\isobj{\Gamma}{O}{\eqclass{S}{O_{1}}{O_{2}}} \\ \isobj{\Gamma}{O'}{\eqclass{S}{O_{1}}{O_{2}}}}
    {\eqobj{\Gamma}{O}{O'}{\eqclass{S}{O_{1}}{O_{2}}}}

  \end{mathpar}

  \caption{Equality Judgments}
  \label{fig:lf-eq}
\end{figure}

\section{Two Type Theories}

The benefit of a logical framework is that it permits the concise specification of type
theories as a signature.

\subsection{G\"odel's T}

\newcommand{\tpsort}{\const{tp}}
\newcommand{\nattp}{\const{nat}}
\newcommand{\arrtp}{\const{arr}}
\newcommand{\arrof}[2]{\appobj{\appobj{\arrtp}{#1}}{#2}}

\newcommand{\elfam}{\const{el}}
\newcommand{\elof}[1]{\appobj{\elfam}{#1}}
\newcommand{\elofp}[1]{\appobj{\elfam}{\parens{#1}}}

\newcommand{\zerocon}{\const{zero}}
\newcommand{\succcon}{\const{succ}}
\newcommand{\succof}[1]{\appobj{\succcon}{#1}}
\newcommand{\reccon}{\const{rec}}
\newcommand{\recof}[4]{\appobj{\appobj{\appobj{\appobj{\reccon}{#1}}{#2}}{#3}}{#4}}

\newcommand{\lamcon}{\const{lam}}
\newcommand{\lamof}[3]{\appobj{\appobj{\appobj{\lamcon}{#1}}{#2}}{#3}}
\newcommand{\appcon}{\const{app}}
\newcommand{\appof}[4]{\appobj{\appobj{\appobj{\appobj{\appcon}{#1}}{#2}}{#3}}{#4}}

\newcommand{\natbetaz}{\const{nat-$\beta$-z}}
\newcommand{\natbetas}{\const{nat-$\beta$-s}}
\newcommand{\arrbeta}{\const{arr-$\beta$}}
\newcommand{\arreta}{\const{arr-$\eta$}}
\newcommand{\pibeta}{\const{pi-$\beta$}}
\newcommand{\pieta}{\const{pi-$\eta$}}

The signature $\Sigma_{T}$ defining G\"{o}del's System T is given in Figure~\ref{fig:t-sig}.
The specifications are verbose, but an implementation would eliminate much of the
redundancy.

\begin{figure}
  
  \begin{align*}
    \tpsort
    & : \sortclass \\
    \elfam 
    & : \arrclass{\tpsort}{\sortclass} \\[1ex]
    \nattp
    & : \tpsort \\
    \arrtp
    & : \arrclass{\tpsort}{\arrclass{\tpsort}{\tpsort}} \\[1ex]
    \zerocon
    & : \elof{\nattp} \\
    \succcon
    & : \arrclass {\elof{\nattp}}{\elof{\nattp}} \\
    \reccon
    &  :
      \piclass{\tpsort}{A}{
      \arrclass{\elof{A}}{
      \arrclass{\parens{\arrclass{\elof{\nattp}}{\arrclass{\elof{A}}{\elof{A}}}}}{
      \arrclass{\elof{\nattp}}{\elof{A}}}}} \\[1ex] 
    \natbetaz
    & :
      \piclass{\tpsort}{A}{
      \piclass{\elof{A}}{b}{
      \piclass{\arrclass{\elof{\nattp}}{\arrclass{\elof{A}}{\elof{A}}}}{s}{
      }}} \\
    & \quad
      \eqclass{\elof{A}}{\recof{A}{b}{s}{\zerocon}}{b} \\
    \natbetas
    & :
      \piclass{\tpsort}{A}{
      \piclass{\elof{A}}{b}{
      \piclass{\arrclass{\elof{\nattp}}{\arrclass{\elof{A}}{\elof{A}}}}{s}{
      }}} \\
    & \quad
      \piclass{\elof{\nattp}}{n}{
      \eqclass{\elof{A}}
      {\recof{A}{b}{s}{\parens{\succof{n}}}}
      {\appobj{\appobj{s}{n}}{\parens{\recof{A}{b}{s}{n}}}}} \\[1ex]
    \lamcon
    & :
      \piclass{\tpsort}{A_{1}}{
      \piclass{\tpsort}{A_{2}}{
      \arrclass{\parens{\arrclass{\elof{A_{1}}}{\elof{A_{2}}}}}{
      \elfam{\parens{\arrof{A_{1}}{A_{2}}}}}}} \\
    \appcon
    & :
      \piclass{\tpsort}{A_{1}}{
      \piclass{\tpsort}{A_{2}}{
      \arrclass{\elofp{\arrof{A_{1}}{A_{2}}}}{
      \arrclass{\elof{A_{1}}}{\elof{A_{2}}}}}} \\[1ex]
    \arrbeta
    & :
      \piclass{\tpsort}{A_{1}}{
      \piclass{\tpsort}{A_{2}}{
      \piclass{\arrclass{\elof{A_{1}}}{\elof{A_{2}}}}{F}{
      \piclass{\elof{A_{1}}}{M_1}{
      }}}} \\
    & \quad
      \eqclass{\elof{A_{2}}}
      {\appof{A_{1}}{A_{2}}{\parens{\lamof{A_{1}}{A_{2}}{F}}}{M_1}}
      {\appobj{F}{M_1}} \\
    \arreta
    & :
      \piclass{\tpsort}{A_{1}}{
      \piclass{\tpsort}{A_{2}}{
      \piclass{\elofp{\arrof{A_{1}}{A_{2}}}}{M}{
      }}} \\
    & \quad
      \eqclass{\elofp{\arrof{A_{1}}{A_{2}}}}
      {M}
      {\lamof{A_{1}}{A_{2}}{\parens{\lamobj{\elof{A_{1}}}{x}{\appof{A_{1}}{A_{2}}{M}{x}}}}}
  \end{align*}

  \caption{Signature of G\"odel's T}
  \label{fig:t-sig}
\end{figure}

\subsection{Dependent T}

The essence of dependent typing is to generalize from types to families of types indexed
by types: if $A:\tpsort$, then $\arrclass{\elofp{A}}{\tpsort}$ is the class of $A$-indexed
families of types.  As a starting point, a reformulation of G\"odel's T in the dependent
setting is given in Figure~\ref{fig:t-sig-dep}.  The main changes are to extend simple
function types to dependent function types, and to generalize the recursor to eliminate
into a $\nattp$-indexed family of result types.  The fully explicit notation is
burdensome, but can be abbreviated in an implementation to permit inference of omitted
parameters and arguments.

\newcommand{\pitp}{\const{pi}}
\newcommand{\piof}[2]{\appobj{\appobj{\pitp}{#1}}{#2}}

\newcommand{\sigtp}{\const{sig}}
\newcommand{\sigof}[2]{\appobj{\appobj{\sigtp}{#1}}{#2}}

\newcommand{\paircon}{\const{pair}}
\newcommand{\pairof}[4]{\appobj{\appobj{\appobj{\appobj{\paircon}{#1}}{#2}}{#3}}{#4}}
\newcommand{\fstcon}{\const{fst}}
\newcommand{\fstof}[3]{\appobj{\appobj{\appobj{\fstcon}{#1}}{#2}}{#3}}
\newcommand{\sndcon}{\const{snd}}
\newcommand{\sndof}[3]{\appobj{\appobj{\appobj{\sndcon}{#1}}{#2}}{#3}}
\newcommand{\splitcon}{\const{split}}
\newcommand{\splitof}[6]{\appobj{\appobj{\appobj{\appobj{\appobj{\appobj{\splitcon}{#1}}{#2}}{#3}}{#4}}{#5}}{#6}}

\begin{exercise}
  Formulate the dependent sum\footnote{aka product, alas} type, in two forms.  The
  formation and introduction rules are the same for the two variants,
  $\sigof{A_{1}}{\parens{\lamobj{\elof{A_{1}}}{x}{A_{2}}}}$ and
  $\pairof{A_{1}}{A_{2}}{M_{1}}{M_{2}}$.  The negative variant takes as elimination forms
  the projections, $\fstof{A_{1}}{A_{2}}{M}$ and $\sndof{A_{1}}{A_{2}}{M}$, with
  appropriate $\beta\eta$ equivalences.  The positive variant takes as elimination form
  splitting, $\splitof{A_{1}}{A_{2}}{B}{M_{1}}{M_{2}}{M}$, with appropriate $\beta\eta$
  equivalences.
\end{exercise}

\begin{figure}
  
  \begin{align*}
    \tpsort
    & : \sortclass \\
    \elfam 
    & : \arrclass{\tpsort}{\sortclass} \\[1ex]
    \nattp
    & : \tpsort \\
    \pitp
    & : \piclass{\tpsort}{A}{\arrclass{\parens{\arrclass{\elof{A}}{\tpsort}}}{\tpsort}} \\[1ex]
    \zerocon
    & : \elof{\nattp} \\
    \succcon
    & : \arrclass {\elof{\nattp}}{\elof{\nattp}} \\
    \reccon
    &  :
      \piclass{\arrclass{\elof{\nattp}}{\tpsort}}{A}{
      \piclass{\elofp{\appobj{A}{\zerocon}}}{b}{
      \piclass{\piclass{\elof{\nattp}}{m}{\arrclass{\elofp{\appobj{A}{m}}}{\elofp{\appobj{A}{\parens{\succof{m}}}}}}}{s}{
      }}} \\
    & \quad
      \piclass{\elof{\nattp}}{n}{\elofp{\appobj{A}{n}}} \\[1ex] 
    \natbetaz
    & :
      \piclass{\arrclass{\elof{\nattp}}{\tpsort}}{A}{
      \piclass{\elofp{\appobj{A}{\zerocon}}}{b}{
      \piclass{\piclass{\elof{\nattp}}{m}{\arrclass{\elofp{\appobj{A}{m}}}{\elofp{\appobj{A}{\parens{\succof{m}}}}}}}{s}{
      }}} \\
    & \quad
      \eqclass{\elofp{\appobj{A}{\zerocon}}}{\recof{A}{b}{s}{\zerocon}}{b} \\
    \natbetas
    & :
      \piclass{\arrclass{\elof{\nattp}}{\tpsort}}{A}{
      \piclass{\elofp{\appobj{A}{\zerocon}}}{b}{
      \piclass{\piclass{\elof{\nattp}}{m}{\arrclass{\elofp{\appobj{A}{m}}}{\elofp{\appobj{A}{\parens{\succof{m}}}}}}}{s}{
      }}} \\
    & \quad
      \piclass{\elof{\nattp}}{n}{
      \eqclass{\elofp{\appobj{A}{\parens{\succof{n}}}}}
      {\recof{A}{b}{s}{\parens{\succof{n}}}}
      {\appobj{\appobj{s}{n}}{\parens{\recof{A}{b}{s}{n}}}}} \\[1ex]
    \lamcon
    & :
      \piclass{\tpsort}{A_{1}}{
      \piclass{\arrclass{\elof{A_{1}}}{\tpsort}}{A_{2}}{
      \arrclass{\parens{\piclass{\elof{A_{1}}}{x}{\elofp{\appobj{A_{2}}{x}}}}}{
      \elfam{\parens{\piof{A_{1}}{A_{2}}}}}}} \\
    \appcon
    & :
      \piclass{\tpsort}{A_{1}}{
      \piclass{\arrclass{\elof{A_{1}}}{\tpsort}}{A_{2}}{
      \arrclass{\elofp{\piof{A_{1}}{A_{2}}}}{
      \piclass{\elof{A_{1}}}{x}{\elofp{\appobj{A_{2}}{x}}}}}} \\[1ex]
    \pibeta
    & :
      \piclass{\tpsort}{A_{1}}{
      \piclass{\arrclass{\elof{A_{1}}}{\tpsort}}{A_{2}}{
      \piclass{\piclass{\elof{A_{1}}}{x}{\elofp{\appobj{A_{2}}{x}}}}{F}{
      \piclass{\elof{A_{1}}}{M_{1}}{
      }}}} \\
    & \quad
      \eqclass{\elofp{\appobj{A_{2}}{M_{1}}}}
      {\appof{A_{1}}{A_{2}}{\parens{\lamof{A_{1}}{A_{2}}{F}}}{M_{1}}}
      {\appobj{F}{M_{1}}} \\
    \pieta
    & :
      \piclass{\tpsort}{A_{1}}{
      \piclass{\arrclass{\elof{A_{1}}}{\tpsort}}{A_{2}}{
      \piclass{\elofp{\piof{A_{1}}{A_{2}}}}{M}{
      }}} \\
    & \quad
      \eqclass{\elofp{\piof{A_{1}}{A_{2}}}}
      {M}
      {\lamof{A_{1}}{A_{2}}{\parens{\lamobj{\elof{A_{1}}}{x}{\appof{A_{1}}{A_{2}}{M}{x}}}}}
  \end{align*}

  \caption{Signature of Dependent G\"odel's T}
  \label{fig:t-sig-dep}
\end{figure}

\subsection{Equality and Identity Types}

Dependent types become interesting only when there are families of types, the principal
examples of which are the extensional and intensional equality types.  Their formulations
are given in Figures~\ref{fig:dep-eq} and~\ref{fig:dep-id}.  The elimination rule for the
extensional equality type is the corresponding equality class whose elimination principle
derives the corresponding equality judgment.\footnote{It is sometimes said that equality
  reflection cannot be formulated in a logical framework, and is therefore suspect.  But
  whether this is so depends on the choice of framework; it is certainly not problematic
  here.}  It follows from this that equality at function type is extensional.  The unicity
rule for equality types states that any two objects of the same equality class are
judgmentally equal; that is, equality classes are ``at most true'' in that the evidence is
immaterial beyond its existence.  The intensional identity type has the same formation and
introduction rules, but has a different elimination rule expressing that the identity type
is the least reflexive relation on the elements of a type.  It is said to be intensional
because it does not validate function extensionality.

\newcommand{\eqtp}{\const{eq}}
\newcommand{\eqof}[3]{\appobj{\appobj{\appobj{\eqtp}{#1}}{#2}}{#3}}
\newcommand{\selfcon}{\const{self}}
\newcommand{\selfof}[2]{\appobj{\appobj{\selfcon}{#1}}{#2}}
\newcommand{\eqrefcon}{\const{eqref}}
\newcommand{\equnicon}{\const{equni}}

\begin{figure}
  
  \begin{align*}
    \eqtp & : \piclass{\tpsort}{A}{\arrclass{\elof{A}}{\arrclass{\elof{A}}{\tpsort}}} \\
    \selfcon & : \piclass{\tpsort}{A}{\piclass{\elof{A}}{M}{\elofp{\eqof{A}{M}{M}}}} \\
    \eqrefcon  & :
         \piclass{\tpsort}{A}{\piclass{\elof{A}}{M_{1},M_{2}}{\arrclass{\elofp{\eqof{A}{M_{1}}{M_{2}}}}{\eqclass{\elof{A}}{M_{1}}{M_{2}}}}}  \\
    \equnicon & :
          \piclass{\tpsort}{A}{\piclass{\elof{A}}{M_{1},M_{2}}{\piclass{\elofp{\eqof{A}{M_{1}}{M_{2}}}}{M,M'}{\eqclass{\elofp{\eqof{A}{M_{1}}{M_{2}}}}{M}{M'}}}}
  \end{align*}

  \caption{Dependent Equality Type}
  \label{fig:dep-eq}
\end{figure}

\newcommand{\idtp}{\const{id}}
\newcommand{\idof}[3]{\appobj{\appobj{\appobj{\idtp}{#1}}{#2}}{#3}}
\newcommand{\reflcon}{\const{refl}}
\newcommand{\reflof}[2]{\appobj{\appobj{\reflcon}{#1}}{#2}}
\newcommand{\jcon}{\const{j}}
\newcommand{\jof}[6]{\appobj{\appobj{\appobj{\appobj{\appobj{\appobj{\jcon}{#1}}{#2}}{#3}}{#4}}{#5}}{#6}}
\newcommand{\idbetacon}{\const{id-$\beta$}}

\begin{figure}
  
  \begin{align*}
    \idtp & : \piclass{\tpsort}{A}{\arrclass{\elof{A}}{\arrclass{\elof{A}}{\tpsort}}} \\
    \reflcon & : \piclass{\tpsort}{A}{\piclass{\elof{A}}{M}{\elofp{\idof{A}{M}{M}}}} \\
    \jcon & :
            \piclass{\tpsort}{A}{
            \piclass{\piclass{\elof{A}}{m_{1}}{\piclass{\elof{A}}{m_{2}}{\arrclass{\elofp{\idof{A}{m_{1}}{m_{2}}}}{\tpsort}}}}{B}{
            }} \\
           & \quad
             \piclass{\piclass{\elof{A}}{x}{\elofp{\appobj{\appobj{\appobj{B}{x}}{x}}{\parens{\reflof{A}{x}}}}}}{r}{
             \piclass{\elof{A}}{m}{\piclass{\elof{A}}{m'}{
             }}} \\
             & \quad
               \piclass{\elofp{\idof{A}{m}{m'}}}{p}{\elofp{\appobj{\appobj{\appobj{B}{m}}{m'}}{p}}} \\
    \idbetacon & :
            \piclass{\tpsort}{A}{
            \piclass{\piclass{\elof{A}}{m_{1}}{\piclass{\elof{A}}{m_{2}}{\arrclass{\elofp{\idof{A}{m_{1}}{m_{2}}}}{\tpsort}}}}{B}{
            }} \\
           & \quad
             \piclass{\piclass{\elof{A}}{x}{\elofp{\appobj{\appobj{\appobj{B}{x}}{x}}{\parens{\reflof{A}{x}}}}}}{r}{
             \piclass{\elof{A}}{m}{
             }} \\
           & \quad
             \eqclass
             {\elofp{\appobj{\appobj{\appobj{B}{m}}{m}}{\parens{\reflof{A}{m}}}}}
             {\jof{A}{B}{r}{m}{m}{\parens{\reflof{A}{m}}}}
             {\appobj{r}{m}}
  \end{align*}

  \caption{Dependent Identity Type}
  \label{fig:dep-id}
\end{figure}

\subsection{Tarskian Universes}

\newcommand{\unitp}{\const{u}}
\newcommand{\unisub}[1]{\appobj{\unitp}{#1}}
\newcommand{\exttp}{\const{ext}}
\newcommand{\extsub}[2]{\appobj{\appobj{\exttp}{#1}}{#2}}
\newcommand{\cumcon}{\const{$\uparrow$}}
\newcommand{\cumsub}[2]{\appobj{\appobj{\cumcon}{#1}}{#2}}
\newcommand{\natcon}{\overline{\const{nat}}}
\newcommand{\picon}{\overline{\const{pi}}}
\newcommand{\piconof}[2]{\appobj{\appobj{\picon}{#1}}{#2}}
\newcommand{\eqcon}{\overline{\const{eq}}}
\newcommand{\eqconof}[3]{\appobj{\appobj{\appobj{\eqcon}{#1}}{#2}}{#3}}
\newcommand{\unicon}{\overline{\const{u}}}
\newcommand{\uniconsub}[1]{\appobj{\unicon}{#1}}
\newcommand{\extnat}{\const{ext-nat}}
\newcommand{\extcum}{\const{ext-cum}}
\newcommand{\extuni}{\const{ext-uni}}
\newcommand{\extpi}{\const{ext-pi}}
\newcommand{\exteq}{\const{ext-eq}}

To add a cumulative hierarchy of universes requires that \LF{} be extended with a class of
natural numbers, written $\lvlclass$, with elements $\zeroobj$ and $\succobj{i}$ for
$i:\lvlclass$.  A Tarskian formulation of universes is given in Figure~\ref{fig:univ}.
Each $\unisub{i}$ is a \emph{universe} whose elements $a$ are \emph{type codes} whose
\emph{extension} as a type is $\extsub{i}{a}$.  Cumulativity is expressed by sending $a$
in $\unisub{i}$ to $\cumsub{i}{a}$ in $\unisub{\parens{\succobj{i}}}$.  The extension of
each of the type codes is defined by equations suggested by the notation.

\begin{figure}[tp]

  \begin{align*}
    \unitp  & : \arrclass{\lvlclass}{\tpsort} \\
    \exttp    & : \piclass{\lvlclass}{i}{\arrclass{\elofp{\unisub{i}}}{\tpsort}} \\
    \cumcon & : \piclass{\lvlclass}{i}{\arrclass{\elofp{\unisub{i}}}{\elofp{\unisub{\parens{\succobj{i}}}}}} \\[1ex]
    \unicon & : \piclass{\lvlclass}{i}{\elofp{\unisub{\parens{\succobj{i}}}}} \\
    \natcon  & : \elofp{\unisub{\zeroobj}} \\
    \picon    & :
                \piclass{\lvlclass}{i}{
                \piclass{\elofp{\unisub{i}}}{a_{1}}{
                \piclass{\arrclass{\elofp{\extsub{i}{a_{1}}}}{\elofp{\unisub{i}}}}{a_{2}}{\elofp{\unisub{i}}}}} \\
    \eqcon  & :
              \piclass{\lvlclass}{i}{
              \piclass{\elofp{\unisub{i}}}{a}{
              \arrclass{\elofp{\extsub{i}{a}}}{\arrclass{\elofp{\extsub{i}{a}}}{\elofp{\unisub{i}}}}}} \\[1ex]
    \extuni & :
              \piclass{\lvlclass}{i}{
              \eqclass{\tpsort}{\extsub{\parens{\succobj{i}}}{\parens{\uniconsub{i}}}}{\unisub{i}}} \\
    \extnat & :
              \eqclass{\tpsort}{\extsub{\zeroobj}{\natcon}}{\nattp} \\
    \extcum & :
              \piclass{\lvlclass}{i}{
              \piclass{\elofp{\unisub{i}}}{a}{
              \eqclass{\tpsort}{\extsub{\parens{\succobj{i}}}{\parens{\cumsub{i}{a}}}}{\extsub{i}{a}}}} \\
    \extpi & :
             \piclass{\lvlclass}{i}{
             \piclass{\elofp{\unisub{i}}}{a_{1}}{
             \piclass{\arrclass{\elofp{\extsub{i}{a_{1}}}}{\elofp{\unisub{i}}}}{a_{2}}{
             }}} \\
            & \quad
              \eqclass
              {\tpsort}
              {\extsub{i}{\parens{\piconof{a_{1}}{a_{2}}}}}
              {\piof{\parens{\extsub{i}{a_{1}}}}{\parens{\lamobj{\elofp{\extsub{i}{a_{1}}}}{x}{\extsub{i}{\parens{\appobj{a_{2}}{x}}}}}}} \\
    \exteq & :
             \piclass{\lvlclass}{i}{
             \piclass{\elofp{\unisub{i}}}{a}{
             \piclass{\elofp{\extsub{i}{a}}}{m_{1}}{
             \piclass{\elofp{\extsub{i}{a}}}{m_{2}}{
             }}}} \\
            & \quad
              \eqclass
              {\tpsort}
              {\extsub{i}{\parens{\eqconof{a}{m_{1}}{m_{2}}}}}
              {\eqof{\parens{\extsub{i}{a}}}{m_{1}}{m_{2}}}
  \end{align*}

  \caption{Tarskian Universes}
  \label{fig:univ}
\end{figure}

\bibliographystyle{plainnat}
\bibliography{notes,martinlof}

\end{document}